\documentclass[12 pt,onesite, openany]{article} 
\usepackage{amsmath,amsxtra,amssymb,latexsym, amscd,amsfonts, indentfirst}
\usepackage[mathscr]{eucal}

\begin{document}
\parskip 3pt
\baselineskip 16pt
\pagenumbering {arabic}
\setcounter{page}{1}
\date{}

\centerline{\bf ON TWO RESULTS OF  MIXED  MULTIPLICITIES } 
 
  \vskip 0.2cm
\centerline {\small Le Van Dinh and Duong Quoc Viet} 

 \centerline {\small Department of Mathematics, Hanoi University of Education} 

 \centerline{\small 136 Xuan Thuy Street, Hanoi, Vietnam} 

 \centerline{\small  E-mail: duongquocviet@fmail.vnn.vn}
  \vskip 0.2cm
 
\noindent This paper shows that the main result of Trung-Verma in 2007 [TV] only is an immediate consequence of an improvement version of     [Theorem 3.4, Vi1]  in 2000.     

Let $(A,\mathfrak{m})$ be a Noetherian local ring of Krull dimension $d=\dim A>0$ with infinite residue field $k=A/\mathfrak{m}.$ Let $J$ be an $\mathfrak{m}$-primary ideal and $I_1,\ldots,I_s$ ideals in $A$ such that their product $I=I_1\cdots I_s$ is non-nilpotent. It is known that the Bhattacharya function 
$\ell_A\Big(\frac{J^nI_1^{n_1}\cdots I_s^{n_s}}{J^{n+1}I_1^{n_1}\cdots I_s^{n_s}}\Big)$
is a polynomial of degree $q-1$ for all sufficiently large $n,n_1,\ldots,n_s$, where $q=\dim A/0:I^\infty$ [Proposition 3.1, Vi1]. If the terms of total degree $q-1$ in this polynomial have the form

\centerline{$\sum_{k_0+k_1+\cdots+k_s=q-1}\frac{1}{k_0!k_1!\cdots k_s!}e(J^{[k_0+1]},I_1^{[k_1]},\ldots,I_s^{[k_s]})n^{k_0}n_1^{k_1}\cdots n_s^{k_s},$}
\noindent then $e(J^{[k_0+1]},I_1^{[k_1]},\ldots,I_s^{[k_s]})$ 
are non-negative integers and are called {\it mixed multiplicities} of the ideals $J,I_1,\ldots,I_s.$

The positivity and the relationship between mixed multiplicities and Hilbert-Samuel multiplicities have attracted much attention. 

Using different sequences, one transmuted  mixed multiplicities into  Hilbert-Samuel multiplicities, for instance: in the case of $\mathfrak{m}$-primary ideals, Risler-Teissier in 1973 [Te] by superficial sequences and Rees in 1984 [Re] by joint reductions; the case of arbitrary ideals, Viet in 2000 [Vi1] by (FC)-sequences and Trung-Verma in 2007 [TV] by $(\varepsilon_1,\ldots,\varepsilon_m)$-superficial sequences.

\vskip 0.2cm
\noindent {\bf Definition 1} [see Definition, Vi1]{\bf.}   
A element $x \in A$ is called an {\it (FC)-element} of $A$ with respect to $(I_1,\ldots, I_s)$ 
if there exists $i \in \{1, 2, \ldots, s\}$  such that $x \in I_i$ and
\begin{list}{}{ \setlength{\leftmargin}{1.8cm}\setlength{\labelwidth}{1.3cm} }\item [\rm (FC1): ] $(x)\cap I_1^{n_1}\cdots I_i^{n_i}\cdots I_s^{n_s} 
= xI_1^{n_1}\cdots I_i^{n_i-1}\cdots I_s^{n_s}$
for all large $n_1, \ldots, n_s$. 
\item [\rm (FC2): ] $x$ is a filter-regular element with respect to $I,$ i.e., $0:x\subseteq 0:I^\infty.$
 \item [\rm (FC3): ] $\dim A/[(x):I^\infty]=\dim A/0:I^\infty-1.$
\end{list} 
We call $x$ a {\it weak-(FC)-element} with respect to $(I_1,\ldots, I_s)$ if $x$ satisfies conditions (FC1) and (FC2).
Let $x_1, \ldots, x_t$ be a sequence in $A$. For each $i = 0, 1, \ldots, t - 1 $, set $A_i = A/(x_1, \ldots, x_{i})S$, $\bar{I}_j = I_j[A/(x_1, \ldots, x_{i})]$, $\bar{x}_{i + 1}$ the image of $x_{i + 1}$ in $A_i$. Then 
$x_1, \ldots, x_t$ is called an {\it(FC)-sequence (respectively, a weak-(FC)-sequence)} of $A$ with respect to $(I_1,\ldots, I_s)$ if $\bar{x}_{i + 1}$ is an (FC)-element (respectively, a weak-(FC)-element)  of $A_i$ with respect to $(\bar{I}_1,\ldots, \bar{I}_s)$ for all $i = 1, \ldots, t - 1$.
\footnotetext{
{\phantom{Key}\it Key words and phrases}: Mixed multiplicity, (FC)-sequence, superficial sequence.}
\vskip 0.2 cm
\noindent {\bf Remark 2.}   
Set $A^*=\dfrac{A}{0:I^\infty}$, $J^* = JA^*$,  ${I_i}^* = I_iA^*$ for all $i= 1,\ldots,s,$ $x \in I_i$ satisfies the condition (FC2) and $x^*$ the image of $x$ in $A^*.$ Then the condition (i) of (Definition in Sect. 3, [Vi1])  is $$(x^*)\cap {I_1^*}^{n_1}\cdots {I_i^*}^{n_i}\cdots {I_s^*}^{n_s} 
= x^*{I_1^*}^{n_1}\cdots {I_i^*}^{n_i-1}\cdots {I_s^*}^{n_s}$$
for all $n_i \geq  n'_i$ and all non-negative integers $n_1, \ldots, n_{i - 1}, n_{i + 1}, \ldots, n_s.$ Since $x$ is an $I$-filter-regular element,  $x^*$ is non-zero-divisor in $A^*.$ Hence        
      $$\frac{x^*{{J^*}^{n_0}}{I_1^*}^{n_1}\cdots {I_i^*}^{n_i-1}\cdots {I_s^*}^{n_s}}
            {x^* {J^*}^{n_0 + 1}{I_1^*}^{n_1}\cdots {I_i^*}^{n_i-1}\cdots {I_s^*}^{n_s}} \cong  
\frac{{J^*}^{n_0}{I_1^*}^{n_1}\cdots {I_i^*}^{n_i-1}\cdots {I_s^*}^{n_s}}
            {{J^*}^{n_0 + 1}{I_1^*}^{n_1}\cdots {I_i^*}^{n_i-1}\cdots {I_s^*}^{n_s}}.$$ Using this property, [Vi1] showed [Proposition 3.3,Vi1]. But in fact, by $x$ satisfies the condition (FC2),       
$$\lambda_x : I_1^{n_1}\cdots I_i^{n_i}\cdots I_s^{n_s} \longrightarrow xI_1^{n_1}\cdots I_i^{n_i}\cdots I_s^{n_s},\;\; y \mapsto xy$$
 is surjective and 
$\text{ker}\lambda_x = (0 : x)\cap I_1^{n_1}\cdots I_i^{n_i}\cdots I_s^{n_s}\subseteq (0 :I^{\infty})\cap I_1^{n_1}\cdots I_i^{n_i}\cdots I_s^{n_s}  = 0$ for all large $n_1, \ldots, n_s$ by Artin-Rees lemma . Therefore, 
$$I_1^{n_1}\cdots I_i^{n_i}\cdots I_s^{n_s} \cong  xI_1^{n_1}\cdots I_i^{n_i}\cdots I_s^{n_s}$$ for all large $n_1, \ldots, n_s.$   This follows that
$$\frac{x\frak J^{n_0}{I_1}^{n_1}\cdots {I_i}^{n_i-1}\cdots {I_s}^{n_s}}
            {x \frak J^{n_0 + 1}{I_1}^{n_1}\cdots {I_i}^{n_i-1}\cdots {I_s}^{n_s}} \cong  
\frac{\frak J^{n_0}{I_1}^{n_1}\cdots {I_i}^{n_i-1}\cdots {I_s}^{n_s}}
            {\frak J^{n_0 + 1}{I_1}^{n_1}\cdots {I_i}^{n_i-1}\cdots {I_s}^{n_s}}$$
for all large $n_1, \ldots, n_s$ and for any ideal $\frak J$ of $A.$ This  proved  that [Proposition 3.3,Vi1] and hence the results of [Vi1] and [Vi2] are still  true for the (FC)-sequences that is defined as in Definition 1.

In this context, Theorem 3.4 in [Vi1] is stated as follows: 
 \vskip 0.2cm
\noindent {\bf Theorem 3} [Theorem 3.4, Vi1]{\bf.} {\it Let $(A, \frak{m})$ denote a Noetherian local ring with maximal ideal 
$\mathfrak{m}$, infinite residue $k = A/\mathfrak{m},$
  and  an ideal  $\frak m$-primary  $J$, and $I_1,\ldots, I_s$   ideals of $A$ such that $I = I_1\cdots I_s$ is non nilpotent. Then the following statements hold.
\begin{enumerate}
\item[\rm (i)] $e(J^{[k_0+1]}, I_1^{[k_1]},\ldots, I_s^{[k_s]}, A) \not= 0$ if and only if there exists an (FC)-sequence 
$x_1,  \ldots, x_t$ $(t = k_1 + \cdots + k_s)$  with respect to $(J, I_1,\ldots, I_s)$ 
consisting of $k_1$ elements of $I_1$, ..., $k_s$ elements of $I_s.$ 
\item[\rm (ii)] Suppose that $e(J^{[k_0+1]}, I_1^{[k_1]},\ldots, I_s^{[k_s]}, A) \not= 0$ and $x_1,  \ldots, x_t$ 
$(t = k_1 + \cdots + k_s)$ is an (FC)-sequence  with respect to $(J, I_1,\ldots, I_s)$ consisting of $k_1$ elements of $I_1$, ..., $k_s$ elements of $I_s$. Set $\bar{A} = A/(x_1, \ldots, x_t):I^\infty$. Then
$$e(J^{[k_0+1]}, I_1^{[k_1]},\ldots, I_s^{[k_s]}, A) = e_A(J,\bar{A}).$$
\end{enumerate}}

Note that Theorem 3 is an immediate cosequence of [Theorem 3.3, VTh] and  the filtration version of Theorem 3 is proved also in [DV].
\vskip 0.2cm
\noindent {\bf Definition 4} [Sect.1, TV]{\bf .}
Set 
$T=\bigoplus_{n_1,\ldots,n_s\geqslant0}\frac{I_1^{n_1}\cdots I_s^{n_s}}{I_1^{n_1+1}\cdots I_s^{n_s+1}}.$
Let $\varepsilon$ be an index with $1\leqslant \varepsilon\leqslant s.$ An element $x\in A$ is an {\it $\varepsilon$-superficial element} for $I_1,\ldots,I_s$ if $x\in I_\varepsilon$ and the image $x^*$ of $x$ in $I_\varepsilon/ I_1\cdots I_{\varepsilon-1}I_\varepsilon^2I_{\varepsilon+1}\cdots I_s$ is a filter-regular element in $T$, i.e., $(0:_Tx^*)_{(n_1,\ldots,n_s)}=0$ for $n_1,\ldots,n_s\gg0.$ Let $\varepsilon_1,\ldots,\varepsilon_m$ be a non-decreasing sequence of indices with $1\leqslant \varepsilon_i\leqslant s.$ A sequence $x_1,\ldots,x_m$ is an {\it $(\varepsilon_1,\ldots,\varepsilon_m)$-superficial sequence} for $I_1,\ldots,I_s$ if for $i=1,\ldots,m$, $\bar{x}_i$ is an $\varepsilon_i$-superficial element for $\bar{I}_1,\ldots,\bar{I}_s$, where $\bar{x}_i,\bar{I}_1,\ldots,\bar{I}_s$ are the images of $x_i,I_1,\ldots,I_s$ in $A/(x_1,\ldots,x_{i-1}).$ 

\vskip 0.2cm
\noindent {\bf Theorem 5} [Theorem 1.4, TV]{\bf.} {\it 
  Set $q=\dim A/0:I^\infty$. Let $k_0,k_1,\ldots,k_s$ be non-negative integers such that $k_0+k_1+\cdots+k_s=q-1.$ Assume that $\varepsilon_1,\ldots,\varepsilon_m$ $(m=k_1+\cdots+k_s)$ is a non-decreasing sequence of indices consisting of $k_1$ numbers $1,\ldots,$ $k_s$ numbers $s$. Let $Q$ be any ideal generated by an $(\varepsilon_1,\ldots,\varepsilon_m)$-superficial sequence for $J,I_1,\ldots,I_s$. Then $e(J^{[k_0+1]},I_1^{[k_1]},\ldots,I_s^{[k_s]})\ne 0$ if and only if $\dim A/Q:I^\infty=k_0+1.$ In this case,
$e(J^{[k_0+1]},I_1^{[k_1]},\ldots,I_s^{[k_s]}) =e\big(J,A/Q:I^\infty\big).$ }

 Then the relationship between  $(\varepsilon_1,\ldots,\varepsilon_m)$-superficial sequences and weak-(FC)-sequences is given in [DV] by the following proposition.
\vskip 0.2cm
\noindent {\bf Proposition 6} [Proposition 4.3, DV]{\bf.}
Let $I_1,\ldots,I_s$ be ideals in $A.$  Let
$x \in A$ be an $\varepsilon$-superficial element for $I_1,\ldots,I_s.$ Then $x$ is a weak-(FC)-element with respect to $(I_1,\ldots, I_s).$  

\vskip 0.2cm
\noindent 
{\bf Proof:}
 Assume that $x$ is an $\varepsilon$-superficial element for $I_1,\ldots,I_s$. Without loss of generality, we may assume that $\varepsilon=1.$ Then
\begin{equation}\label{eq3}
\big(I_1^{n_1+2}I_2^{n_2+1}\cdots I_s^{n_s+1}:x\big) \cap I_1^{n_1}\cdots I_s^{n_s} = I_1^{n_1+1}I_2^{n_2+1}\cdots I_s^{n_s+1}
\end{equation}
for $n_1,\ldots,n_s\gg0.$ (1) implies
\begin{equation}\label{eq6}
\begin{split}
\big(I_1^{n_1+2}I_2^{n_2+1}\cdots I_s^{n_s+1}:x\big) \cap I_1^{n_1}I_2^{n_2+1}\cdots I_s^{n_s+1}
 = I_1^{n_1+1}I_2^{n_2+1}\cdots I_s^{n_s+1}
\end{split}
\end{equation}
for $n_1,\ldots,n_s\gg0$. We prove by induction on $k\geqslant2$ that
\begin{equation}\label{eq7}
\begin{split}
\big(I_1^{n_1+k}I_2^{n_2+1}\cdots I_s^{n_s+1}:x\big) \cap I_1^{n_1}I_2^{n_2+1}\cdots I_s^{n_s+1}
= I_1^{n_1+k-1}I_2^{n_2+1}\cdots I_s ^{n_s+1}
\end{split}
\end{equation}
for $n_1,\ldots,n_s\gg0.$ The case $k=2$ follows from (\ref{eq6}). Assume now that
\begin{equation*}
\begin{split}
\big(I_1^{n_1+k}I_2^{n_2+1}\cdots I_s^{n_s+1}:x\big) \cap I_1^{n_1}I_2^{n_2+1}\cdots I_s^{n_s+1}
=I_1^{n_1+k-1}I_2^{n_2+1}\cdots I_s^{n_s+1}
\end{split}
\end{equation*}
for $n_1,\ldots,n_s\gg0.$ Then 
$$\begin{aligned}
\big(&I_1^{n_1+k+1}I_2^{n_2+1}\cdots I_s^{n_s+1}:x\big) \cap I_1^{n_1}I_2^{n_2+1}\cdots I_s^{n_s+1}\\
 &= \big(I_1^{n_1+k+1}I_2^{n_2+1}\cdots I_s^{n_s+1}:x\big)
 \cap \big(I_1^{n_1+k}I_2^{n_2+1}\cdots I_s^{n_s+1}:x\big) \cap I_1^{n_1}I_2^{n_2+1}\cdots I_s^{n_s+1}\\
&= \big(I_1^{n_1+k+1}I_2^{n_2+1}\cdots I_s^{n_s+1}:x\big)\cap  I_1^{n_1+k-1}I_2^{n_2+1}\cdots I_s^{n_s+1}\\
&=I_1^{n_1+k}I_2^{n_2+1}\cdots I_s^{n_s+1}\end{aligned}$$ 
for $n_1,\ldots,n_s\gg0.$ The last equality is derived from (\ref{eq6}). Hence the induction is complete and we get (\ref{eq7}). It follows that for $n_1,\ldots,n_s\gg0$,
$$\begin{aligned}
(0:x)& \cap I_1^{n_1}I_2^{n_2+1}\cdots I_s^{n_s+1}\\ 
&=\Big(\bigcap_{k\geqslant2}I_1^{n_1+k}I_2^{n_2+1}\cdots I_s^{n_s+1}:x\Big) \cap I_1^{n_1}I_2^{n_2+1}\cdots I_s^{n_s+1}\\
&=\Big(\bigcap_{k\geqslant2}\big(I_1^{n_1+k}I_2^{n_2+1}\cdots I_s^{n_s+1}:x\big)\Big) \cap I_1^{n_1}I_2^{n_2+1}\cdots I_s^{n_s+1}\\
&=\bigcap_{k\geqslant2}\Big(\big(I_1^{n_1+k}I_2^{n_2+1}\cdots I_s^{n_s+1}:x\big) \cap I_1^{n_1}I_2^{n_2+1}\cdots I_s^{n_s+1}\Big)\\
&=\bigcap_{k\geqslant2}I_1^{n_1+k-1}I_2^{n_2+1}\cdots I_s^{n_s+1}=0,
\end{aligned}$$
that is, $(0:x) \cap I^n = 0$ for $n \gg0,$ here $I = I_1\cdots I_s.$  Hence $0:x\subseteq 0:I^\infty.$  So  $x$ is satisfies condition (FC2). Now we need  to prove that
$I_1^{n_1}\cdots I_s^{n_s}\cap(x)=xI_1^{n_1-1}I_2^{n_2}\cdots I_s^{n_s}$
for $n_1,\ldots,n_s\gg0$. But this has from the proof of [Lemma 1.3, TV]. Hence $x$ is a weak-(FC)-element with respect to $(I_1,\ldots,I_s)$.

\vskip 0.2cm
\noindent {\bf Remark 7.} Return to Theorem 5, assume that $Q=(x_1,\ldots,x_m)$, where $x_1,\ldots,x_m$ is an $(\varepsilon_1,\ldots,\varepsilon_m)$-superficial sequence for $J, I_1,\ldots,I_s$. As $x_1,\ldots,x_m$ is a weak-(FC)-sequence with respect to  $(J, I_1,\ldots,I_s)$ by Proposition 6. Hence 

\centerline{$\dim A/Q:I^\infty\leqslant q-m=k_0+1$}
\noindent with equality if and only if $x_1,\ldots,x_m$ is an (FC)-sequence by [Proposition 3.1(ii), Vi2].    
This fact proved  that Theorem 3 covers Theorem 5 that is the main result of Trung and Verma in [TV].

\end{document}